\magnification\magstep1
\input amstex

\input amsppt.sty
\vsize22truecm
\hoffset1.3truecm
\voffset2.4truecm
\parindent20pt
\font\got=eufm10
\def\t{\text}

\def\aj{{\alpha_j}}
\def\cic{C^{(m)}}
\def\cim{C^{m-1}}
\def\cis{C^{m-2}}
\def\tic{T^{(m)}}
\def\goto{\text{\got S}}
\def\ov{\overline}
\def\rt{\Cal R}
\def\cu{\Cal C}
\NoBlackBoxes

\topmatter
\title How many latin rectangles are there?\endtitle
\author Aurelio de Gennaro\endauthor

\abstract Until now the problem of counting Latin rectangles 
$m\times n$ has been 
solved with an explicit formula for $m=2,3$ and 4 only.
In the present paper an explicit formula is provided for the
calculation of the number of Latin rectangles for any
order $m$. The results attained up to now become
particular cases of this new formula. Furthermore, putting $m=n $,
the number of Latin squares of
order $n$ can also be obtained in an
explicit form.\endabstract\endtopmatter

\head 0. Introduction\endhead

A Latin rectangle $m\times n$  is a matrix with $n$ rows and $n$
columns the elements of which are chosen in 
$[n]=\{1,\dots,n\}$ so that two elements are never the same,
neither  on the same row nor on the same column. It is said that
such a Latin rectangle has order $m$. From the definition it
follows that $m\le n$ and that each row of a Latin rectangle is a
permutation of $[n]$.

Furthermore it is clear that it is always possible to standardize the
first row making it the same as the permutation $12\dots
n$. In such a case we say that we are dealing with a ``reduced Latin
rectangle''.

If we call the number of Latin rectangles $m\times n$
with $L(m,n)$ and the number of reduced Latin
rectangles with  the same dimension with $K(m,n)$, it is clear that
$L(m,n)=n!K(m,n)$.

The problem of counting Latin rectangles has engaged several
generations of  mathematicians, but the results reached up to now, 
as we will see later,  are limited to certain special cases.

With this paper we intend to finally supply an explicit formula for
the calculation of $K(m,n)$ for any value of the order   $m$.

The partial results attained  up to now will result special
cases of such a formula. 

The result obviously also allows  the
calculation of the number  of Latin squares of order $n$,
which represent  the special case of $m\times n$ Latin
rectangles in which $m$ takes  on its maximum admissible value 
$n$.

\head 1. A brief survey of results\endhead

A Latin rectangle consists of $m$ permutations of $[n]$ which, taken
two by two, don't have fixed points. It is from this point of view
that the problem was initially studied by Montmort, Euler and
Lucas.

It seems that the solution in the simplest case $m=2$,
 known as ``derangement problem'', can be found going back
to Montmort 1713 [6]. It consists of the number $D_n$ of the
permutations of $[n]$ without fixed points given by:
$$
D_n=\sum^n_0 {}\!_{{}_k}  (-1)^k{n!\over k!}\tag1.1 $$ 
and equivalent to  $K(2,n)$.

In 1891 Lucas expounded the famous ``m\'enage problem'' which
consists of counting the ways of arranging $n$ couples at a round
table so that men and women alternate and no husband and wife
are adjacent to one another. The problem, examined since 1878 by
Tate, was also studied by Cayley and Muir, however no satisfactory
results were reached.

The solution to m\'enage problem is equivalent to the
enumeration of all the permutations of $[n]$ which are
discordant with both the permutations $12\dots n$ and $23\dots
n1$.
The generalisation of the above mentioned problem --- known as
``the cyclical Touchard problem of index $m$'' or ``problem of
$m$-discordant permutations'' --- proposes the counting of all
the permutations $\sigma$ of $[n]$ so that: $\forall\, i\in [n],\
\sigma(i)-i\not\equiv 0,1,\dots,m-1$ (mod $n$).

It is also said that this problem is equivalent to the enumeration of
the ``very reduced Latin rectangles'' of order $m+1$, the number
of which is indicated with  $V(m,n)$, that have the first $m$
permutations $\sigma_j$, with $j\in[m]$, in the canonical form: 
$\sigma_j(i)=i+j-1$ (mod $n$).

The solution in the simplest  case $m=2$ was found by Touchard in
1934 [15] and consists of  Touchard's famous numbers 
$U_n$, equivalent to $V(2,n)$, expressed by: 
$$
U_n=\sum^n_0{}\!_{{}_k}(-1)^k{2n\over 2n-k}{2n-k\choose k}
(n-k)!\ .\tag1.2 
$$
For the subsequent case $V(3,n)$ recursive algorithms have been
obtained by Riordan [12] and by Yamamoto [18]. However an
explicit formula was only provided in 1967 by Moser [7]. In the case
$V(4,n)$, the biggest yet dealt with, there is only one
recursive result by Whitehead in 1979 [16] and possibly an
explicit formula by Nechvatal [8] also in 1979.

However let us return to the more general and more complex
problem of the calculation of  $K(m,n)$, which represents the aim
of this paper.

The first attempts at the calculation of  $K(3,n)$ go back to 
Jacob and to Kerawala who, in 1941 [5], found a recursive
formula. The following tidy explicit formula for $K(3,n)$ is, on the
other hand, attributed to Yamamoto (see [1] and [10]): 
$$
K(3,n)=n!\sum_{a+b+c=n}(-1)^b\,2^c\,{a!\over c!}{3a+b+2\choose
b}\ .\tag1.3
$$
Furthermore, in 1944, Riordan obtained an expression of $K(3,n)$
in terms of Touchard's numbers $U_n$ and subsequently, in 
1946 [11], the well known formula:
$$
K(3,n)=\sum_0^{\left[{n\over 2}\right]}{}\!_{{}_k}{n\choose
k}D_k\,D_{n-k}\,U_{n-2k}\tag1.4
$$
(with $U_0=1$) which expresses $K(3,n)$ in terms of  $D_k$
and $U_k$. 

It is necessary to say that until now, in this line of
research, no other progress has been achieved since, for $m>3$, it 
has not been possible to obtain $K(m,n)$ in terms of $K(i,n)$ and
$V(j,n)$ with $i,j<m$.

The case of $K(4,n)$, which is the most complex yet to be dealt with
successfully, was only solved with an explicit formula
in 1979, this was achieved independently by Nechvatal [8] and by
Athreya, Pranesachar and  Singhi [1].

Subsequently, in 1980, Pranesachar [10] and Nechvatal [9],
by different means, found a way to express
$K(m,n)$ for any value of $m$ by means of the M\"obius
function of the lattice of partitions of a set. The limit of these
research works is that they take the calculation of $K(m,n)$
back to the enumeration of other combinatorial objects, such as
the partitions of an integer, for which no explicit
formulas are known, and thus don't allow an explicit formula for
$K(m,n)$ to be obtained. A further tidy result of this type was
achieved by Gessel in 1987 [3]. 

Until now, then, no explicit formula is known which permits
the calculation of $K(m,n)$ whatever the value of $m$.

We would like to conclude this section remembering that
another interesting line of research tried to get asymptotic
expressions of  $K(m,n)$. The first significant paper of this kind is
attributed to  Erd\"os and Kaplansky [2] in 1946, subsequent
results were obtained by Yamamoto in 1951 [17] and by Stein in
1978 [14].  

Finally, in recent years,  Godsil and McKay [4] achieved an
asymptotic valuation of $V(m,n)$.

\head 2. Notation and preliminaries\endhead

Very useful concepts in the study of permutations without fixed
points are those of board and of rook polynomial.

A board is a nonempty subset of $\Bbb
P\times\Bbb P$ ($\Bbb P=$ set of positive integers), the elements
of board are called squares. Considering a board  $C$ it is usually
indicated with $r_k(C)$ the number of different ways of placing  $k$
non-attacking rooks on it.
The rook polynomial of $C$ in the symbolic variable $x$, that we'll
call $R(C)$, is given by $\sum_k
r_k(C)x^k$. 

If we write  $\goto_n$ for the set of all permutations
of $[n]$, then every $\sigma \in \goto_n$ can be thought of as
a board, called the ``graph'' of $\sigma$, the squares of which
are the couples $\left(i,\sigma(i)\right)$ $\forall\, i\in[n]$. It is
furthermore obvious that $m$ permutations, two by two without
fixed points, make up a board of  $m\cdot n$ squares, if in such a
board the first permutation is the identical one $\sigma(i)=1$, it will
hereon  be indicated with $C^{(m)}$.

We furthermore state, to have a greater number of symbols,
that a number put up to the right of a symbol doesn't
denote a raising to a power of the same, but it acts as a new
symbol (therefore $a^2$ it isn't the square of  $a$). When we want
to indicate $a$ raised to $m$, we write $(a)^m$.

We will call $C^g$ the board, included in $C^{(m)}$, formed by the
$g$-th permutation of  $C^{(m)}\ \left(g\in [m]\right)$. Since, as
has been mentioned previously, to speak about  $C^{(m)}$
it is the same as to speak about a reduced Latin rectangle of order
$m$,  we will refer often to  $C^g$  as the  $g$-th line of $C^{(m)}$
(which is  not to be confused with the $g$-th row or column of 
$C^{(m)}$ like a board that are different things). We can even say
that a subset of  $C^g$ has ``grade'' $g$.  

At this point let us remember a classic result which joins the
rook polynomials to the permutations without fixed points.We
consider the permutations as boards and, taking the board 
$B\subseteq [n]\times[n]$, we indicate with $N_s(B)$ the number
of permutations of $[n]$ which have exactly $s$ squares in common
with $B$. And so giving us the following tidy relation:
$$
N_s(B)=\sum^n_s{}\!_{{}_k}(-1)^{k-s}{k\choose s}(n-k)!
r_k(B)\tag2.1
$$
for the proof of this see [13] chapter 2.3.

Now let us introduce some conventions in the use of symbols. The
sets are always indicated with capital letters and the number of
the elements which make them up with  the corresponding lower
case letter  (therefore $a=|A|$). $C(A)$ will be the complementary
of the set $A$ in the universe set. 

$T^{(m)}$ will represent a generic system of independent rooks ---
that is which don't attack each other --- put on  $\cic$ and $T^g$,
with $g\in [m]$, will be the part of the system contained on the
$g$-th line of $\cic\ \left(T^g=\tic\cap C^g\right)$ and thus it will
be:  $\tic=T^1\cup T^2\cup\cdots\cup T^m$.

Furthermore, if $A\subseteq\cic $, we will say that $\Cal R(A)$ is
the projection for rows of $A$ in $\cic$, consisting of all the squares
of  $\cic$ which have any square of  $A$ in their own row
(obviously $A\subseteq \Cal R(A)$). We will also say that $\Cal
R_g(A)$ is the projection for  rows of  $A$ on the $g$-th line of
$\cic$ and we will put $\Cal R_g(A)=\Cal R(A)\cap C^g$, with the
consequence that: $\Cal R(A)=\Cal R_1(A)\cup\cdots\cup \Cal
R_m(A)$.  Similarly, speaking of columns instead of rows, we can
define  $\Cal C(A)$ and $\Cal C_g(A)$. 

Finally we define the set  $\Cal I(A)=\Cal
R(A)\cap\Cal C(A)$ as the ``impression'' of $A$ and the set $\Cal
O(A)=\Cal R(A)\cup \Cal C(A)$ as the ``shadow'' of $A$. 
Thus
$\Cal O(\tic)$ will be the set of all the squares of  $\cic$ subject to
the attack of any rook of $\tic$ and which therefore can't contain
other independent rooks from those of $\tic$.

To indicate the number of ways in which the set $A$ can
generally be arranged, considering the restrictions which have
been imposed on it, we will write $\pi(A)$. 
So we shall obtain that
$K(m,n)=\pi(\cic)$. 

As it is known, a partition in $k$ blocks of a set $A$
is formed  by a collection of non empty sets $A_i$, with $i\in
[k]$,  two by two disjoint and such
that $\bigcup\limits^k_1{}\!_i\,A_i=A$. We will indicate with  $\prod(A)$
the set of the partitions of $A$; if $\pi\in\prod(A)$ and $\pi$ has $k$
blocks, we say that  $|\pi|=k$; finally
we put $\prod_n=\prod\left([n]\right)$.

Let us also remember that the refinement of two partitions 
$\pi_1$ and $\pi_2$, with $\pi_1,\pi_2\in \prod(A)$, is the partition
of  $A$ consisting of all the non empty intersections of some block of 
$\pi_1$ with some block of $\pi_2$.

If $X=\{x_1,\dots,x_s\}$ is a set of variables, we put, for economy of
space: 
$$
{n\choose X}={n\choose x_1,\dots,x_s}
$$
and furthermore:
$$
\sum X=\sum_{x_i\in X}x_i\ ;\  \prod X=\prod_{x_i\in X}x_i\ ;\
\prod X!=\prod_{x_i\in X}x_i !
$$
and, ranging each $x_i$ within its own domain:
$$
\sum_X=\sum_{x_1}\dots\sum_{x_s}\ .
$$
We will indicate with $(n)_k=n(n-1)\dots(n-k+1)$ the falling
factorial of  $n$ and with $\langle n\rangle_k=n(n+1)\dots(n+k-1)$
the raising factorial of $n$. Now, given that $(n)_n=n!$,
we intend to put $\langle n\rangle_n=n$!` even if in other
literature this symbolism has been used to indicate the
subfactorial of  $n$.

Let us conclude this section with some recalls relative to
the permutations of $[n]$.

If $\sigma \in\goto_n $ and $\sigma (i)=a_i\ \forall\, i\in[n]$, 
we can also say that $\sigma$ corresponds to the word
$a_1\dots a_n$. It is well known, see [13] pages 17 and following,
that $\sigma $ can be shared in an unambiguous way in the
product of disjoint cycles on the elements of  $[n]$  and that it can
have a ``standard representation'', which we will indicate with 
$s_1\dots s_n$, writing a) the elements of each cycle with the
largest element first and  b) arranging the cycles in increasing
order of their largest element. In such a case each cycle will start
with a left--to--right maximum i. e.  with an element $s_i$ so
that $s_i>s_j$ for each $j<i$.

Another way to describe  $\sigma$ can be achieved by indicating
with $b_i$ the number of elements $j$ of its standard
representation on the left of $i$ with $j>i$ and defining 
$I(\sigma)=(b_1,\dots,b_n)$ the ``inversion table'' of  $\sigma$. In
fact it can be easily proved, see [13] Proposition 1.3.9, that
between the $\sigma\in\goto_n$ and the $I(\sigma)$ there is a
bijection and furthermore that: $0\le b_i\le n-i,\ \forall\,i\in[n]$.

\head 3. The associated partitions\endhead

In this section we want to show how the computation of 
$K(m,n)$ can be taken back to the enumeration of  a double
system of partitions of $[n]$. 

The first step in the argument is to reiterate the method of
calculation by means of systems of independent rooks expressed
by the formula (2.1).

Let us suppose that we want to determine the number $K(m,n)$ of
all the possible $\cic$ and that we have already counted all the
possible arrangements of the first $m-1$ lines $C^{(m-1)}$, the
number of ways in which $C^{m}$ can be arranged can be obtained
easily, by means of the formula (2.1), once
$r_k\left(C^{(m-1)}\right)$,  for $k=0,\dots, m-1$, are known.

In fact, if $T^{(m-1)}_{1/m-1}=T^1_{1/m-1}\cup\cdots\cup
 T^{m-1}_{1/m-1}$ is the independent generic system of rooks on
$C^{(m-1)}$ and therefore $t^{(m-1)}_{1/m-1}
=t^1_{1/m-1}+\cdots+t^{m-1}_{1/m-1}=k$, we will have that: 
$$
\pi(C^m)=\sum^n_0\!{}_{t^{(m-1)}_{1/m-1}}
(-1)^{{}^{t^{(m-1)}_{1/m-1}}}\left(n-t^{(m-1)}_{1/m-1}\right)!
\pi\left(T^{(m-1)}_{1/m-1}\right)\ .\tag3.1
$$
Now we are trying to calculate $\pi(C^{m-1})$ with the same
assumptions.
Now
$\pi(C^{m-1})=\pi(T^{m-1}_{1/m-1})\pi\left(C^{m-1}-T^{m-1}_{1/m-1}
\right)$ and, if we consider the generic 
$T^{(m-2)}_{1/m-2}\subseteq
 \Cal I\left(T^{m-1}_{1/m-1}\right)$ and
$T^{(m-2)}_{2/m-2}\subseteq \Cal
I\left(C^{m-1}-T^{m-1}_{1/m-1}\right)$, we will have that: 
$$
\align
& \pi(\cim)=\sum_{t^{(m-2)}_{1/m-2}}
\sum_{t^{(m-2)}_{2/m-2}}
(-1)^{t^{(m-2)}_{1/m-2}}
\left({t^{m-1}_{1/m-1}}-{t^{(m-2)}_{1/m-2}}\right)!\tag3.2\\
&\hskip-1cm
\cdot\pi\left(T^{(m-2)}_{1/m-2}\right)
(-1)^{t^{(m-2)}_{2/m-2}}\left(n-
{t^{m-1}_{1/m-1}}-{t^{(m-2)}_{2/m-2}}\right)!\pi
\left(T^{(m-2)}_{2/m-2}\right)
.\endalign
$$
Repeating the argument for the subsequent line  $\cis$, we see
that: 
$$
\align
&\cis=T^{m-2}_{1/m-2}\cup\left(T^{m-2}_{1/m-1}-
T^{m-2}_{2/m-2}\right)
\cup\biggl(T^{m-2}_{1/m-1}
\cap
T^{m-2}_{2/m-2}\biggr)\tag3.3\\
&\hskip-1cm\cup\left(T^{m-2}_{2/m-2}-T^{m-2}_{1/m-1}\right)\cup\left(
\cis-(T^{m-2}_{1/m-2}\cup
T^{m-2}_{1/m-1}\cup T^{m-2}_{2/m-2})\right)
\endalign
$$
and, calling: $T^{(m-3)}_{1/m-3},\ 
T^{(m-3)}_{2/m-3},\ T^{(m-3)}_{3/m-3},\
T^{(m-3)}_{4/m-3},\ T^{(m-3)}_{5/m-3}$
the generic $T^{(m-3)}$ included in the impression of each of the
five components of the union of which at (3.3), we will be able to
say, referring always to (2.1), that: 
$$
\align
&\pi(\cis)=\pi\left(T^{m-2}_{1/m-2}\right)
\pi\left(T^{m-2}_{1/m-1}
-T^{m-2}_{2/m-2}\right)
 \tag3.4\\
&\cdot\pi\biggl(T^{m-2}_{1/m-1}\cap T^{m-2}_{2/m-2}\biggr)
\pi\left(T^{m-2}_{2/m-2}-T^{m-2}_{1/m-1}\right)
\pi\biggl(\cis-T^{m-2}_{1/m-1}\\
&\cup T^{m-2}_{1/m-2}\cup T^{m-2}_{2/m-2}\biggr)=
\sum_{t^{(m-3)}_{i/m-3}}(-1)^{\sum\limits^5_1\!{}_i
\,t^{(m-3)}_{i/m-3}}
\left(t^{m-2}_{1/m-2}-t^{(m-3)}_{1/m-3}\right)!
\endalign
$$
$$
\align
&\cdot\left(|T^{m-2}_{1/m-1}-T^{m-2}_{2/m-2}|
-t^{(m-3)}_{2/m-3}\right)!
\left(|T^{m-2}_{1/m-1}\cap T^{m-2}_{2/m-2}|
-t^{(m-3)}_{3/m-3}\right)!\\
&\cdot\left(|T^{m-2}_{2/m-2}-
T^{m-2}_{1/m-1}|
-t^{(m-3)}_{4/m-3}\right)!\bigl(|\cis-\left(T^{m-2}_{1/m-1}\cup
T^{m-2}_{1/m-2}\cup T^{m-2}_{2/m-2}\right)|\\
&-t^{(m-3)}_{5/m-3}\bigr)!
\pi\left(T^{(m-3)}_{1/m-3}\right)
\pi\left(T^{(m-3)}_{2/m-3}\right)
\pi\left(T^{(m-3)}_{3/m-3}\right)
\pi\left(T^{(m-3)}_{4/m-3}\right)
\pi\left(T^{(m-3)}_{5/m-3}\right)\ .
\endalign
$$
Continuing in this way in order to count all the possible
arrangements of the line $l$, we should take all the possible subsets
of rooks which are situated on it, considering the partition
refinement of  $C^l$ consisting of all their non empty intersections
$T^l_j$, $j=1,\dots,p_l$, and of the complementary of their union
$T^l_0=C^l-\bigcup\limits^{p_l}_1\!{}_j T^l_j$ and finally choosing  
$p_l+1$
systems of independent rooks $T^{(l-1)}_{j/l-1}$, $j=0,1,\dots,p_l$,
with $T^{(l-1)}_{j/l-1}\subseteq \Cal I(T^l_j)$.

Then we will have:
$$
\pi(C^l)=\prod^{p_l}_0\!{}_{{}_j}\pi(T^l_j)=
\sum_{t^{(l-1)}_{j/l-1}}(-1)^{\sum\limits^{p_l}_0\!{}_{{}_j}
t^{(l-1)}_{j/l-1}} 
\prod^{p_l}_0\!{}_{{}_j}(t^l_j-t^{(l-1)}_{j/l-1})!\pi(T^{(l-1)}_{j/l-1})
\tag3.5
$$
and we can conclude that:
$$
K(m,n)=\prod^m_2\!{}_{{}_l}\,\pi(C^l)
=\sum_{t^{(l-1)}_{j/l-1}}(-1)^{{}^{\sum\limits^m_2\!{}_l\sum
\limits^{p_l}_0\!{}_j
t^{(l-1)}_{j/l-1}}}\prod^m_2\!{}_{{}_l}\prod^{p_l}_0\!{}_{{}_j}\left(t^l_j-
t^{(l-1)}_{j/l-1}\right)!
\pi(T^{(l-1)}_{j/l-1})\ .\tag3.6
$$
The situation, therefore, seems to be somewhat complex, but an
idea which allows us to control it is that of identifying any subset of
$\cic$ by means of its two projections, for rows and for columns,
on the main diagonal $C^1=\{(i,i)\ i\in[n]\}$.

Given a set of grade $g$ $A^g$,  we consider in fact $\Cal
R_1(A^g)$ and $\Cal C_1(A^g)$. Now, if $g=1$, $\Cal R_1(A^1)=\Cal
C_1(A^1)=A^1$ and  projections and set coincide. If, on the
other hand, $g>1$, then the set $A^g$ determines
$\Cal R_1(A^g)$ and $\Cal C_1(A^g)$ in one way only. Viceversa
given $A^g_1$,  $A^g_2\subseteq C^1$ with $|A^g_1|=|A^g_2|=|A^g|=p$,
$A^g$ will be one of the  $p!$ permutations of the square board
$\Cal R\left(A^g_1\right)\cap \Cal C\left(A^g_2\right)$.
Furthermore, if  $B$ is a board of forbidden positions for
$A^g$ and we set $\ov B=B\cap \Cal R\left(A^g_1\right)\cap
\Cal C\left(A^g_2\right)$, we will have that:
$\pi(A^g)=\sum\limits_0^{p}\!{}_k(-1)^k(p-k)!r_k(\ov B)\ .$

In the light of this new approach, a system of independent rooks
$T^{(l)}=T^1\cup\cdots\cup T^l$ determines the subsets
$R^i=\Cal R(T^i)$ and $C^i=\cu(T^i)$, for $i\in [l]$, and:
$R^0=\rt\left(C(T^{(l)})\right)=C\left(\bigcup\limits^l_1\!{}_i
R^i\right)$ and $C^0=\cu\left(
C(T^{(l)})\right)=C\left(\bigcup\limits^l_1\!{}_i C^i\right)$.  Thus it
characterizes two partitions, each one with $l+1$ blocks, the
$\{R^i\}$ and the  $\{C^i\}$, with $i=0,\dots,l$, the first for the set of
the rows and the second for that of the columns. 

Now, if we return to the computation performed with (3.5) of all
the possible arrangements of $C^{(l)}$, we will have that
$T^{(l-1)}_{j/l-1}$ characterizes a partition in  $l$ blocks
$R^i_{j/l-1}=\Cal R(T^i_{j/l-1})$, with $i=0,\dots, l-1$, of the block
$\rt (T^l_j)$ of the partition $\{\rt(T^l_j)\},\ j=0,\dots, p_l$, in
$p_l+1$  blocks of the set of the rows.

Consequently, if we put together all the elements of grade $i$  of
the various partitions $R^i_{j/l-1} $ and  we set
$R^i_l=\bigcup\limits^{p_l}_0\!{}_j R^i_{j/l-1}$, we will  obtain that
$\{R^i_l\}$, with $i=0,\dots, l-1$, is a partition in  $l$ blocks
of the set of the rows of our board. Similarly we can construct the
partition $\{C^i_l\}$, with $i=0,\dots, l-1$, of the set of columns.
Repeating this for each line from the  $m$-th to the second, we will
eventually have  $m-1$ couples of partitions of the set of the rows
and of that of the columns: $\{R^i_j\}$ and  $\{C_j^i\}$, with
$j=2,\dots, m$ and $i=0,1,\dots j-1$.

Intersecting each of these partitions with $C^1$ we obtain as many
partitions of  $C^1=\{(i,i)\}$ with $i\in[n]$. In such partitions the
projections for rows and for columns of the sets of grade 1, which
belong to  $C^1$, obviously coincide. 

Calculating the number of all
the possible arrangements of these  $2(m-1)$  partitions, 
respecting the condition that the projections of the sets of grade 1
must coincide,  is equivalent, from what has been said, to
calculating the product
$\prod\limits^m_2\!{}_{{}_l}\prod\limits^{p_l}_0\!{}_{{}_j}\pi
\left(T^{(l-1)}_{j/l-1}\right)$ which appears in  (3.6).

To do this it is natural to consider the partition refinement of the
$m-1$ $\{R^i_j\cap C^1\}$ partitions and that of the  
$m-1$ $\{C^i_j\cap C^1\}$ partitions. Putting by analogy
$R^0_1=C^0_1=C^1$, the blocks of the refinement partitions will be 
given from: $R_{\alpha_m,\dots,\alpha_1}=R^{\alpha_m}_m\cap
R^{\alpha_{m-1}}_{m-1}\cap\cdots\cap R^{\alpha_1}_1$, with $0\le
\alpha_j\le j-1\ \forall\,j\in [m]$, and from
$C_{\beta_m,\dots,\beta_1}=C^{\beta_m}_m\cap
C^{\beta_{m-1}}_{m-1}\cap\cdots \cap C_1^{\beta_1}$ with similar
limitations on the indices $\beta_j$.

We will then have:
$$
\prod^m_2\!{}_{{}_l}\prod^{p_l}_0\!{}_{{}_j}\pi
(T^{(l-1)}_{j/l-1})=\prod_0^{j-1}\!{}_{{}_{\alpha_j}}
\pi(R_{\alpha_m,\dots,\alpha_1})
\prod_0^{j-1}\!{}_{{}_{\beta_j}}\pi(C_{\beta_m,\dots,\beta_1})\
.\tag3.7 $$
We will also say that the pair of partitions 
$\{R_{\alpha_m,\dots,\alpha_1}\}$ and  $\{C_{\beta_m,\dots,\beta_1}\}$ of
$[n]$ is ``associated'' with the collection of systems of
independent rooks $T^{(l-1)}_{j/l-1}$, with $l=2,\dots, m$ and $
j=0,\dots, p_l$.

\head 4. The blocks of the associated partitions\endhead

Before being able to develop the calculation of (3.6) using
(3.7), we must examine closely the meaning of the indices 
$\alpha_m,\dots, \alpha_1$ and $\beta_m,\dots, \beta_1$ which
respectively mark the blocks $R_{\alpha_m,\dots,\alpha_1}$ and
$C_{\beta_m,\dots,\beta_1}$ .

In order to do this, we first of all define the concept of ``covering''.
Taking the index  $\alpha_i$, with $\alpha_i>0$, we say that
$\alpha_i$ ``covers'' $\alpha_{\alpha_i}$ and we write 
$\alpha_i\vdash \alpha_{\alpha_i}$ or
$\kappa(\alpha_i)=\alpha_{\alpha_i}$. \medskip

\noindent >From the definition it immediately follows that:
\item{a)} if $\alpha_i=0$ it doesn't cover any other index, since
$\alpha_0$ doesn't exist; 
\item{b)} if $\alpha_i\vdash \alpha_s$
then $s<i$.

\noindent Thus, if $\alpha_s>0 $ , applying to it more times
the function of covering $\kappa$, you will always arrive at an
index of value 0.

On the contrary, taking an index  $\alpha_l$ of value 0, we can
consider all the indices $\alpha_s$ which have the property
$\alpha_s\vdash \alpha_l$ (that is $\kappa^{-1}(\alpha_l)$).  Repeating
this procedure more times, we get all the indices which, with a
finite number of applications of the function $\kappa$, finish in 
$\alpha_l$. 

If we suppose that $\alpha_h=0$  and we put 
$\kappa^0(\alpha_s)=\alpha_s$,
we will be able to define  $Z_h=\{\alpha_j\,|\, \alpha_h=0\ \text{\rm
and } \exists\, k\in\Bbb N\ \text{\rm so that 
}\kappa^k(\alpha_j)=\alpha_h\}$ which we will call the ``component h''
of the indices $\alpha_m,\dots, \alpha_1$.
Furthermore, as $\alpha_1=0$, we will have that $Z_1\ne
\emptyset$.

In this way, we obtain a partition of the set of indices 
$\{\alpha_m,\dots,\alpha_1\}$ in blocks made up from  $Z_h$,
with $h\in [m]$.

We will say that such a subdivision represents the structure  of the
indices $\alpha_m,\dots,\alpha_1$ and we will write that
$\sigma(\alpha_m,\dots,\alpha_1)=Z_1\cup Z_{z_2}\cup\dots\cup
Z_{z_a}$ with $1<z_2<\cdots<z_a\le m$. We will also write 
$\sigma_l(\alpha_m,\dots,\alpha_1)=Z_l$ and
$\zeta(\alpha_m,\dots,\alpha_1)=$number of the $\alpha_j$ which are
equal to zero. It is clear that, if $Z=\{\alpha_{i_1},\dots,\alpha_{i_s}\}$ 
is a component, then $\alpha_{i_s}=0$.

Moreover we will put, to be brief
$R_{\alpha_j}=R_{\alpha_m,\dots,\alpha_1}$,
$C_{\beta_j}=C_{\beta_m,\dots,\beta_1}$, 
$\alpha_j=\alpha_m,\dots,\alpha_1$ and
$\beta_j=\beta_m,\dots,\beta_1$.

It is necessary to pay attention to the fact that  $Z_h$ is not only a
subset  $I\subseteq [m]$, but a subset of the indices $\alpha_j$,
for $j\in I$, each with its own value.

The following result allows to count the number of the 
$R_{\alpha_j}$ at the base of the structure of their indices.
\proclaim{4.1 Proposition}  Let $\sigma(\alpha_m,\dots,
\alpha_1)=Z_1\cup Z$ and $Z=Z_{z_2}\cup\cdots\cup Z_{z_a}$,
then:
\item{a)} if we suppose $Z$ to be variable, the number of the
possible sets of indices will be: ${m-1\choose z}z!$. 
\item{b)} if, on the other hand, we keep $Z$ constant, then the
possible $\alpha_j$ will be $(m-1-z)!$.\endproclaim

\noindent In fact, to determine $Z$ we will have,  first of all, have
to choose the $z$ places of its indices in the set  $\{2,\dots,m\} $
and this can be done in  ${m-1\choose z}$ ways.
Furthermore, if the selected indices are
$\alpha_{j_1},\dots,\alpha_{j_z}$ (with $j_1>\cdots>j_z$), it can be
seen that $\alpha_{j_z}$ must be equal to 0, $\alpha_{j_{z-1}}$
can assume the values  0 and $\alpha_{j_z}$ and so on. Therefore
the last index has only one possible value, the penultimate two
values etc., thus all the possible ways to attribute a value to 
$\alpha_{j_1},\dots,\alpha_{j_z}$ are $1\cdot2\cdot\dots \cdot z=z!$.  
And this proves a).

If, on the other hand,  $Z$ is fixed, the places of the indices of $Z_1$
are also fixed. Now the last index of $Z_1$ on the right $\alpha_1$
can only have the value 0, the penultimate only 1 and thus, for an
argument identical to the previous, the number  of possible values
of the $m-z$ indices is equal to  $1\cdot1\cdot
2\cdot\dots\cdot(m-1-z)=(m-1-z)!$.
And thus b) too  is proved. 

It is also possible to calculate the number of possible $\aj$ , in
terms of the data of singular components with the following result,
which we shall just state.

\proclaim {4.2 Proposition} If
$\sigma(\alpha_m,\dots,\alpha_1)=Z_1\cup
Z_{z_2}\cup\cdots\cup Z_{z_a}$, then $ |Z_1| $ and $|Z_{z_i}|$, with
$i=2,\dots,a$, constitute a partition of the integer $m$ in which
the number of parts equal to  $s$ will be $\lambda_s$. The number 
of possible  $\alpha_j$ with this structure will therefore be the
same as: 
$$ m!\over 1^{\lambda_1}2^{\lambda_2}\dots
m^{\lambda_m}\  \lambda_1!\lambda_2!\dots\lambda_m!\tag4.1
$$\endproclaim\bigskip

\noindent So, if $\alpha_s\vdash \alpha_l$, we have that
$R_{\alpha_j}$ is a subset of
$\rt_1\left(\bigcup_x\,T^{\alpha_l}_{x/l-1}\right)$, with $x$ that
ranges in a subset of  $\{0,1,\dots,p_l\}$, and thus it lies in the
projection for rows on the first line of  a set of rooks of grade 
$\alpha_l$ included in the  impression of a set of rooks of grade
$\alpha_s=l$. If instead $\alpha_s$  doesn't cover $\alpha_l$, then
the set of rooks of grade  $\alpha_l$  is not included in the
impression of the set of rooks of grade $\alpha_s$ and so the
number of elements in their intersection varies according to the
variation of the projection for rows or for columns.

From this, it follows that if, using the symbolism of  section 
3, we take
$\rt_1(T^l_j)-\rt_1(\bigcup\limits^{l-1}_1\!{}_i T^i_{j/l-1})$, 
we see that it
will be composed of the union of all the $R_{\alpha_j}$ with the
same  $Z_l$ component. Viceversa, if we fix the 
$Z_l$  component and make the other $\alpha_j$  vary in all
possible ways, we obtain a collection of sets  $R_{\alpha_j}$ the
union of which will be equal to 
$\rt_1(T^l_j)-\rt_1(T^{(l-1)}_{j/l-1})$ for some  $j$. Furthermore, if 
$l=1$, since $T^{(l-1)}_{j/l-1}$ doesn't exist, the union of all the 
$R_{\alpha_j}$ with the same  $Z_1$ will be given by $\rt_1(T^1_j)$
for some $j$.

Naturally the same  argument  is true for the sets 
$C_{\beta_j}$ and the components of the indices $\beta_j$.

\head 5. The enumeration of Latin rectangles\endhead

Now let us try, applying the contents of the previous section,
to give an explicit form to (3.6) in terms of the
data of the two associated partitions $\{R_{\alpha_j}\}$  and
$\{C_{\beta_j}\}$, and that is in terms of the sets of variables
$R=\{r_{\alpha_j}\}$ and $C=\{c_{\beta_j}\}$.

First of all, we observe that, if
$l>1$,
$t^l_j-t^{(l-1)}_{j/l-1}=|T^l_j|-|T^{(l-1)}_{j/l-1}|=|\rt_1(T^l_j)|-|\rt_1
(T^{(l-1)}_{j/l-1})|=|\rt_1(T^l_j)-\rt_1(T^{(l-1)}_{j/l-1})|$ but,
following what was said before,
$\rt_1(T^l_j)-\rt_1(T^{(l-1)}_{j/l-1})$ is formed, in such a case,
from the union of all the sets
$R_{\alpha_j}$ with the same $Z_l$ component and viceversa.

Therefore, if we put, $\forall \,l\in [m]$,
$Q(Z_l)=\{r_{\alpha_j}\,\vert\,\sigma_l(\alpha_j)=Z_l\},\ 
\tilde Q(Z_l)=\{c_{\beta_j}\,|\,r_{\beta_j}\in Q(Z_l)\}$
and $q(Z_l)=\sum\,Q(Z_l)$, thus we have that, if $l>1$, $j\in
\{0,\dots,p_l\}$ exists  so that:
$$
q(Z_l)=t^l_j-t^{(l-1)}_{j/l-1}\tag5.1
$$
and viceversa.

Now, we will compute the product
$\prod\limits^m_2\!{}_{{}_l}\prod\limits^{p_l}_0\!{}_{{}_j}
\,\pi(T^{(l-1)}_{j/l-1})$
that, for (3.7), is the same as
$\prod\limits^{j-1}_0\!{}_{{}_{\alpha_j}}
\prod\limits^{j-1}_0\!{}_{{}_{\beta_j}}\,
\pi(R_{\alpha_j})\, \pi(C_{\beta_j})$ .

The first partition  $\{R_\aj\}$ can be  chosen in a completely
arbitrary way and thus the number of its possible arrangements
is given by the multinomial coefficient ${n\choose R}$. The
second partition is, on the other hand, subject to some restrictions.

First of all, for  $l>1$, $T^{(l-1)}_{j/l-1}\subseteq \Cal
I(T^l_j)$  and so $|\cu_1(T^l_j)-\cu_1(T^{(l-1)}_{j/l-1})|=
|\rt_1(T^l_j)-\rt_1(T^{(l-1)}_{j/l-1})|=q(Z_l)$ 
and since, following the same reasoning as we have already done,
$|\cu_1 (T^l_j)-\cu_1(T^{(l-1)}_{j/l-1})|=\sum\,\tilde Q(Z_l)$
we have that:
$$\sum\,\tilde Q(Z_l)=q(Z_l)\ .\tag5.2
$$
Furthermore, taking a generic set of rooks of grade  $l$  it is clear
that $|\cu_1(T^l)|=|\rt_1(T^l)|$. 

Now, if  $l=1$, then  $T^l\subseteq C_1$ and even
$\cu_1(T^l)=\rt_1(T)^l$, but, for the reasons stated in section 4, 
$\rt_1(T)^1$ is made up of the union of all the $R_{\alpha_j}$
with the same $Z_1$,  and the same argument is valid for 
$\cu_1(T^1)$, therefore:
$$
\bigcup_{\sigma_1(\alpha_j)=Z_1}R_{\alpha_j}=
\bigcup_{\sigma_1(\beta_j)=Z_1}C_{\beta_j}\tag5.3
$$
and:
$$
\sum\,\tilde Q(Z_1)=q(Z_1)\ .\tag5.4
$$
Furthermore, the restrictions (5.2) and (5.4) imposed on 
$c_{\beta_j}$ imply that, for  $T^l_j$, with $j=0,\dots,p_l$,
$|\cu_1(T^l_j)|=|\rt_1(T^l_j)|$.

This can be easily proved for complete induction on $l$
considering that, if $l=1$, the result has already been expressed
by (5.4), while, if $l>1$ and we suppose that we have already
proved this $\forall\, j\in[l-1]$, it follows from the consideration
that: $|T^l_j|=\sum\limits^{l-1}_1\!{}_i\,|\,
T^i_{j/l-1}|+|\rt_1(T^l_j)-\rt_1(T^{(l-1)}_{j/l-1})|$.

Therefore there are no other restrictions on  $c_{\beta_j}$, apart
from those expressed by (5.2) and (5.4).

If we now group the $C_{\beta_j}$ sets on the basis of the
value of their component $Z_1$, (5.3) allows us to state that:
$$
\prod_{\beta_j}\pi(C_{\beta_j})=\prod_{Z_1}{q(Z_1)\choose \tilde
Q(Z_1)}\tag5.5
$$
on the condition, however, that the  $C$ variables also respect the
restrictions imposed by (5.2).

Let us finally examine the $t^{(l-1)}_{j/l-1} $ which appear in  (3.6)
as exponents of $-1$.

Now, for  (5.1), if $l>1$, $t^{(l-1)}_{j/l-1}=t^l_j-q(Z_l)$ for any
$Z_l$ and so $\sum\limits^{p_l}_0\!{}_{{}_j}
t_{j/l-1}^{(l-1)}=\sum\limits^{p_l}_0\!{{}_j}t_j^l-\sum_{{}_{Z_l}}
q(Z_l)=n-\sum_{{}_{Z_l}}q(Z_l)$.
Therefore, being  $l=2,\dots,m$, the exponent of  $-1$ in  (3.6) will
be the same as 
$n(m-1)-\sum\limits^m_2\!{}_{{}_l}\sum_{{}_{Z_l}}q(Z_l)$.
Furthermore, since, as we  have already seen, 
$n=\sum\,R=\sum_{Z_1}q(Z_1)$, adding and subtracting  $n$ it can
be expressed  by:
$nm-\sum\limits^m_1\!{}_{{}_l}\sum_{{}_{Z_l}}q(Z_l)$.

Finally, set $W=\{r_{\alpha_j}\,|\,\zeta(\alpha_j)\ \text{\rm
odd}\}$ and considering that every $r_{\alpha_j}$ 
variable  compares $\zeta(\alpha_j)$ times in 
$\sum\limits^m_1\!{}_{{}_l}\sum_{{}_{Z_l}}q(Z_l)$, we will have that 
the exponent
of  $-1$ can be substituted by  $nm-\sum W$, since the even
multiples of  $r_{\alpha_j}$ can obviously be omitted.

Using all these results in  (3.6), we obtain the following remarkable
result: 
\proclaim{5.1 Theorem} 
$$
\multline
K(m,n)=\sum\!{}_{{}_R}\sum_{\sum\tilde
Q(Z_l)=q(Z_l)}\hskip-.8truecm{}_{{}_C}\ (-1)^{{}^{nm+
\sum\limits^m_1\!{}_l\sum_{Z_l}
q(Z_l)}}
{n\choose R}\prod^m_2\!{}_l\prod_{Z_l}q(Z_l)!\\
\cdot\prod_{Z_1}{q(Z_1)\choose \tilde Q(Z_1)}=\sum_{\sum\,
R=n}\!{}_{{}_R}\sum_{\sum\tilde Q(Z_l)=q(Z_l)}
\hskip-.8truecm{}_{{}_C}\ 
(-1)^{{}^{nm+\sum\,W}}{\prod\limits^m_0\!{}_l\prod_{Z_l}q(Z_l)!\over
\prod R!\prod C!}\\
\hskip-2,5truecm=(-1)^{n(m-1)}\sum_{\sum R=n}\!{}_{{}_R}\sum_{\sum\tilde
Q(Z_l)=q(Z_l)}\hskip-.8truecm{}_{{}_C}\ {\prod\limits^m_0\!{}_l\prod_{Z_l}\left(-q(Z_l)
\right)\t{\rm !`}\over
\prod R!\prod C!}\\
=(-1)^{nm}\sum\!{}_{{}_R}\sum{}_{{}_C}\ {n\choose
R}\prod^m_1\!{}_{{}_l}\prod\!{}_{{}_{Z_l}}(-1)^{q(Z_l)}{q(Z_l)\choose
\tilde Q(Z_l)}\cdot\prod_C c_{\beta_j}^{\zeta(\beta_j)-1}\ .
\endmultline\tag5.6
$$
\endproclaim
Where, by analogy with the preceding symbolism, we have set
$Z_0=\emptyset$ since  $\alpha_0$ doesn't exist. Thus
$q(Z_0)=\sum R=n$ since the elements of  $Q(Z_0)$, not being
subject to any restrictions, are all the elements of $R$.

So  (5.6) is an explicit formula for the computation of $K(m,n)$ in
$2m!$ variables $R$ and $C$, while $q(Z_l)$ with $l\in [m]$ and $\sum
W$ are sums of particular subsets of $R$.

This therefore represents the result which  we proposed to achieve
with the present paper.

The $C$ variables, in contrast to the  $R$ variables, are not,
however, between their independent since  they must be subject
to the restrictions $\sum\tilde Q(Z_l)=q(Z_l)$ for $l\in[m]$.

If we want to limit ourselves to considering only independent
variables we can proceed as follows.

For each $Z_l$ component we indicate with  $d(Z_l)$ the 
$c_{\beta_j}$ variable with $\sigma_l(\beta_j)$ $=Z_l$ and all
the indices $\beta_j$ which are different from those of $Z_l$ equal
to zero, and we put $D=\{d(Z_l)\}$.

Now $d(Z_l)=q(Z_l)-\sum\limits_{{}_{c_{\beta_j}\in \tilde
Q(Z_l)-D}}c_{\beta_j}$ and thus the variables of  $D$ can be
obtained from those of  $C-D$.

Furthermore, if  $\sigma_l(\beta_j)=Z_l$, then $c_{\beta_j}\in\tilde
Q(Z_l)$ and so $c_{\beta_j}\le q(Z_l)$. Thus, if we put 
$\mu_{\beta_j}=\min_{{}_{\sigma_l(\beta_j)\ne\emptyset}}
\left(q(\sigma_l(\beta_j))\right)$,
then $\forall\,c_{\beta_j}\in C-D$, we will have that 
$c_{\beta_j}\le\mu_{\beta_j}$ and such a restriction
guarantees that $d(Z_l)\ge 0$.

Using this new symbolism (5.6) can be rewritten like this:
$$
K(m,n)=\sum^n_0\!{}_{{}_{r_{\alpha_j}}}\sum_0^{\mu_{\beta_j}}
\!{}_{{}_{c_{\beta_j}}}
(-1)^{{}^{nm+\sum W}}{n\choose
R}{\prod\limits^m_1\!{}_{{}_l}\prod_{{}_{Z_l}}q(Z_l)!\over
\prod(C-D)!\prod D!}\tag5.7
$$
with  $r_{\alpha_j}\in R$ and $c_{\beta_j}\in C-D$.

Now if  $|Z_l|=s$, for the  Proposition 4.1, the possible $d(Z_l)$ are
${m\choose s}(s-1)!$. Furthermore, if $s=1$, all the $d(Z_l)$
coincide with the $c_{\beta_j}$ which has all the indices at 0 and so,
 in such a case, instead of ${m\choose 1}0!=m$ we only have 
one distinct element and $|D|=\sum\limits^m_1\!{}_s{m\choose
s}(s-1)!-(m-1)$.

Thus, in  (5.7), other than the  $m!$ independent variables $R$,
there are the $m!+m-1-\sum\limits^m_1\!{}_s{m\choose s}(s-1)!$
independent variables $C-D$.

\head 6. Simplifications of the formula\endhead

We have seen that (5.7) needs  
$2m!+m-1-\sum\limits^m_1\!{}_s{m\choose s}(s-1)!$ 
independent variables for the computation of $K(m,n)$.

It is possible, though, to effect two types of elimination among
these parameters which allow us to reduce their number
considerably, even though this fact makes  (5.7)  lose its
symmetry. This is obviously important when we would like to
calculate concretely  $K(m,n)$ for $m$ and $n$
prefixed.

Let us therefore examine the two possible reductions of  the
independent variables $R$ and $C-D$.

\noindent A) We consider  $r_{\alpha_j}$ and $c_{\alpha_j}$ with
$\zeta (\alpha_j)=1$ and so with $\sigma(\alpha_j)=Z_1$. In such an
assumption  $q\left(\sigma(\alpha_j)\right)$ contains a unique
element and so, for (5.4), $c_{\alpha_j}=r_{\alpha_j}$ and, in  (5.7),
$c_{\alpha_j}!$ is simplified with
$q\left(\sigma(\alpha_j)\right)!=r_{\alpha_j}!$.
As far as $r_{\alpha_j}$ is concerned instead, if we put:
$F_0=\{r_{\alpha_j} \,|\,\zeta(\alpha_j)=1\}\ ,\ f_0=\sum F_0,\ov
Q_0=R-F_0$  and $\tilde F_0=\{c_{\alpha_j}|r_{\alpha_j}\in F_0\}$,
we will have that the variables of $F_0$ don't appear in any set
$Q(Z_l)$ with $l>1$ and that: $$
\align
&\sum_R{n\choose R}=\sum_{F_0}\sum_{\ov Q_0}{n\choose
F_0,\ov Q_0}=\sum_{F_0}\sum_{\ov Q_0}{n\choose f_0}{f_0\choose
F_0}{n-f_0\choose \ov Q_0}\tag 6.1\\ 
&\qquad=\sum_{\ov
Q_0}\sum_{f_0}\left((m-1)!\right)^{f_0}{n!\over
f_0!\prod \ov Q_0!}
\endalign
$$
since, for the Proposition 4.1, $|F_0|=(m-1)!$. Furthermore the
$F_0$ appear among the exponents of $-1$ with their total
$f_0$. The  $2(m-1)!$ variables of $F_0$ and of $\tilde F_0$
can therefore be substituted by $f_0$.

\noindent B)  Let us now consider the $r_{\alpha_j}$ and
$c_{\alpha_j}$  with $\sigma(\alpha_j)=Z_h\cup Z_1$ (and so
$\zeta(\alpha_j)=2$) and $\min(z_1,z_h)=1$ and put,
$\forall\,s\in[m]$: $F_s=\{r_{\alpha_j}\,|\,\zeta(\alpha_j)=2\
\text{\rm and }|\sigma_s(\alpha_j)|=1\}$, $f_s=\sum\,F_s$,
$\ov Q(Z_s)=Q(Z_s)-F_s$, $\tilde
F_s=\{c_{\alpha_j}\,|\,r_{\alpha_j}\in F_s\}$ and  $\ov
q(Z_s)=q(Z_s)-f_s$. First we observe that, if
$r_\aj\in F_s$, it doesn't appear among the exponents of $-1$ since 
$\zeta(\aj)$ is even. Now, if $\sigma(\aj)=Z_s\cup Z_v$, $|Z_v|=m-1$
and so $q(Z_v)$ has only one element and, for (5.2) and (5.4),
$c_\aj=r_\aj$. Therefore in  (5.7), if $c_\aj \in \tilde F_s, \  c_\aj!$ is
simplified with $q(Z_v)!=r_\aj!$. Moreover,  $\forall\, s\in[m]$: 
$$
\align
&\sum_{F_s}{q(Z_s)!\over \prod F_s!}=\sum_{F_s}{\left(
f_s+\ov q(Z_s)\right)!\over \prod F_s!}\tag6.2\\
&\hskip-1cm=\sum_{F_s}{\left(
f_s+\ov q(Z_s)\right)!\over f_s!}{f_s\choose
F_s}=\left((m-2)!\right)^{f_s} {\left(
f_s+\ov q(Z_s)\right)!\over f_s!}
\endalign
$$
since $|F_s|=(m-2)!$, and so also the $F_s$ and the $\tilde F_s$
are eliminated and substituted by $f_s$.
We must, however, by careful because, if $m=2$, $F_1$ and $F_2$
are equal.

In conclusion, putting $\ov R=R-\bigcup\limits^m_0\!{}_iF_i$, $\ov
C=\{c_\aj\,|\,r_\aj\in\ov R\}$ and $\ov
D=\{d\left(\sigma_l(\aj)\right)\,|\,c_\aj\in \ov C\}$ we have that
(5.7) transforms itself into:
$$
\align
&K(m,n)=\sum_{f_0}\sum_{f_s}\sum^n_0\!{}_{{}_{\ov
R}}\sum^{\mu_{\beta_j}}_0\!{}_{{}_{\ov C-\ov D}}(-1)^{{}^{nm+\sum
W}}\left((m-1)!\right)^{f_0}\left((m-2)!
\right)^{\sum\limits^m_1\!{}_sf_s}\tag6.3\\
&\quad\cdot{1\over
\prod (\ov C-\ov D)!\prod \ov D!}
{n\choose f_0,\dots, f_m,\ov
R}\prod^m_1\!{}_{{}_s}\left(f_s+\ov
q(Z_s)\right)!\prod^m_1\!{}_{{}_l}\prod^{m-2}_2\!{}_{{}_{|Z_l|}} 
q(Z_l)!\ .
\endalign
$$
It is, however, possible to accomplish a further step to simplify
(6.3). In fact, putting $f=\sum\limits^m_1\!{}_sf_s $, we have that:
$$
\align
&\sum_{f_1}\cdots\sum_{f_m}\prod^m_1\!{}_{{}_s}{(f_s+\ov
q(Z_s))!\over 
f_s!}=\sum_{f_1}\cdots\sum_{f_m}\prod^m_1\!{}_{{}_s}\ov
q(Z_s)!{f_s+\ov q(Z_s)\choose \ov q(Z_s)}\tag6.4\\
&\quad=\prod^m_1\!{}_{{}_s}\ov
q(Z_s)!{f+\sum\limits^m_1\!{}_s\ov q(Z_s)+m-1\choose f}
\endalign
$$
and so (6.3) becomes:
$$
\align
&K(m,n)
=\sum^n_0\!{}_{{}_{f_0}}\sum^n_0
\!{}_{{}_{f}}\sum^n_0\!{}_{{}_{\ov R}}
\sum^{\mu_{\beta_j}}_0
\!{}_{{}_{\ov C-\ov D}}(-1)^{{}^{nm+\sum
W}}\left((m-1)!\right)^{f_0}\tag6.5\\
&\quad\cdot\left((m-2)!\right)^{f}
{n\choose f_0,  
f,\ov R}\,{1\over \prod(\ov C-\ov D)!\prod\ov D!}\\ 
&\quad\cdot{\left(f+\sum\limits^m_1\!{}_s\ov
q(Z_s)+m-1\right)!\over
\left(\sum^m_1\!{}_s\ov
q(Z_s)+m-1\right)!}\prod^m_1\!{}_l
\prod^{m-2}_2\!{}_{{}_{|Z_l|}}q(Z_l)!
\prod^m_1\!{}_{{}_s}\ov
q(Z_s)!\ .
\endalign
$$
From the independent variables $R$ we have so eliminated the 
$(m-1)!$ of $F_0$  and the $(m-2)!$ of each $F_s$, with $s\in [m]$,
and therefore $|\ov R|=m!-(m-1)!-m(m-2)!=m!-(2m-1)(m-2)!$.

The $C$ variables have undergone the same reduction. However it  is
necessary to add $f_0$ and  $f$ and subtract the $\ov D$, which are
as many as the components $Z_l$ with $1<|Z_l|<m-1$, and so equal to
$\sum\limits^{m-2}_2\!{}_{{}_h}{m\choose h}(h-1)!$ plus the $c_\aj$
with all the $\aj$ indices equal to zero  (which is determined by
the $m$ equivalent restrictions $\ov q(Z_s)=\sum\left(\tilde
Q(Z_s)-\tilde F(Z_s)\right)$ with $s\in [m]$), and thus:  
$$
\align
&|\ov
C-\ov
D|=m!-(m-1)!-m(m-2)!-\biggl(\sum^{m-2}_2\!{}_h{m\choose
h}(h-1)!\tag6.6\\
&\quad+1\biggr)
=m!-(2m-1)(m-2)!-1
-\sum^{m-2}_2\!{}_h
{m\choose h}(h-1)!\\
&\quad=m!-1
-\sum^m_2\!{}_h{m\choose h}(h-1)!
\ .\endalign
$$
The independent parameters of (6.5) are therefore all together:
$2m!-(2m-1)(m-2)!+m+1-\sum\limits^m_1\!{}_h{m\choose
h}(h-1)!\ .$

\head 7. The simplest cases\endhead

Let us see what in concrete terms happens calculating
the formulas obtained in sections 5 and 6 for the
first values  of $m=2,3,4$.

\noindent A) $m=2$. $r_\aj$ are of 
$r_{\alpha_2\alpha_1}$ which can therefore assume the values  
$r_{10}$ and $r_{00}$. Furthermore $C-D=\emptyset$,
$c_{00}=r_{00}$ and $c_{10}=r_{10} $ and so, applying (5.7),
we have: 
$$
K(2,n)=\sum^n_0\!{}_{{}_{r_{10}}}\sum^n_0\!{}_{{}_{r_{00}}}
(-1)^{{}^{2n+r_{10}}}r_{00}!{n\choose
r_{10},r_{00}}{r_{00}!r_{10}!\over 
c_{00}!c_{10}!}=\sum_0^n\!{}_{{}_{r_{10}}}(-1)^{r_{10}}{n!\over
r_{10}!}\tag7.1
$$
which, for  (1.1), is equivalent to $D_n$. 

\noindent B) $m=3$.  $r_\aj$ are of
$r_{\alpha_3\alpha_2\alpha_1}$. As seen in section 6, the
$F_0=\{r_{210},r_{110}\},\ F_1=\{r_{200}\},\
F_2=\{r_{100}\},\ F_3=\{r_{010}\}$ and the homologous  $c_\aj$
are eliminated. Furthermore $\ov q(Z_1)=\ov q(Z_2)=\ov
q(Z_3)=r_{000}$ and $c_{000}=r_{000}$ and so, applying (6.5),
we have that: 
$$
\align
&K(3,n)=\sum^n_0\!{}_{f_0}\sum^n_0\!{}_{f}\sum^n_0
\!{}_{r_{000}}
(-1)^{3n-f_0-r_{000}}\,2^{f_0}\,1^{f}\,{n\choose
f_0,f,r_{000}}\tag7.2\\
&\quad\cdot{(r_{000}!)^3(f+3r_{000}+2)!\over
c_{000}!(3r_{000}+2)!}
=\sum_{f_0+f+r_{000}=n}(-1)^{f}\, 2^{f_0}\,{n!r_{000}!\over
f_0!}{3r_{000}+f+2\choose f} \endalign
$$
and we find (1.3) again.

\noindent C) $m=4$. $r_\aj$ are of
$r_{\alpha_4\alpha_3 \alpha_2\alpha_1}$. The $F_0=\{r_{3210},
r_{2210},r_{1210}, r_{1110},r_{2110},r_{3110}\}$, 
$F_1=\{r_{3200},r_{2200}\}$,
$F_2=\{r_{3100},r_{1100}\}$,  
$F_3=\{r_{1010},r_{2010}\}$,
$F_4=\{r_{0110},r_{0210}\}$ and the homologous $c_\aj$ are
eliminated. Furthermore $q(Z'_1)=r_{1000}+r_{1200}$, 
$q(Z''_1)=r_{2100}+r_{0100}$, 
$q(Z'''_1)=r_{0010}+r_{3010}$, 
$q(Z'_2)=r_{0200}+r_{1200},\
q(Z''_2)=r_{2000}+r_{2100},\
q(Z_3)=r_{3000}+r_{3010},$
and $\ov
q(Z_1)=r_{0000}+r_{3000}+r_{2000}+r_{0200},\
\ov
q(Z_2)=r_{0000}+r_{3000}+r_{0100}+r_{1000}$,
$\ov
q(Z_3)=r_{0000}+r_{2000}+r_{1000}+r_{0010},\ 
\ov
q(Z_4)=r_{0000}+r_{0100}+r_{0010}+r_{0200}$.
We besides have that:
$c_{1000}=r_{1000}+r_{1200}-c_{1200}$,
$c_{0100}=r_{0100}+ r_{2100}-c_{2100}$, 
$c_{0010}=r_{0010}+r_{3010}-c_{3010}$, 
$c_{0200}=r_{0200}+
r_{1200}-c_{1200}$,
$c_{2000}=r_{2000}+r_{2100}-c_{2100}$,
$c_{3000}=r_{3000}+ r_{3010}-c_{3010}$,
$c_{0000}=r_{0000}+r_{3000}+r_{2000}+
r_{0200}-c_{3000}-c_{2000}-c_{0200} $,  that: 
$\ov
R=\{r_{1000},r_{1200},r_{2100},r_{0100},r_{0010},
r_{3010},r_{0200},r_{3000},r_{2000},r_{0000}\}$
and that: $\ov C-\ov D=\{c_{1200},c_{2100},c_{3010}\}$
and, applying (6.3), we obtain:
$$
\multline
K(4,n)=\sum_{f_0}
\sum_{f_s}\sum_{\ov R}\sum_{\ov C-\ov
D}
(-1)^{4n-f_0-r_{1000}-r_{0100}-r_{0010}-r_{3000}
-r_{2000}-r_{0200}}\\
6^{f_0}2^{f_1+f_2+f_3+f_4}
{n\choose f_0,f_1,f_2,f_3,f_4,\ov
R}(f_2+r_{0000}+r_{3000}
+r_{0100}+r_{1000})!\\
\cdot(f_3+r_{0000}+r_{2000}
+r_{1000}+r_{0010})!
(f_4+r_{0000}+r_{0100}
+r_{0010}+r_{0200})!\\
\cdot{(f_1+r_{0000}+r_{3000}+r_{2000}
+r_{0200})!\over
c_{0000}!c_{3000}!c_{2000}!c_{0200}!} (r_{0200}+r_{1200})!
(r_{2000}+r_{2100})!\\
\cdot(r_{3000}+r_{3010})!{r_{1000}+r_{1200}\choose 
c_{1000},c_{1200}}
{r_{2100}+r_{0100}\choose 
c_{2100},c_{0100}}
{r_{0010}+r_{3010}\choose 
c_{0010},c_{3010}}\qquad
\endmultline\tag7.3$$
$$
\multline
=\sum_{{f_0}+\sum{f_s}+\sum{\ov R}=n}\,\sum_{\ov C-\ov
D}(-1)^{f_0+r_{1000}+r_{0100}+r_{0010}+
r_{3000}+r_{2000}+r_{0200}}\\
6^{f_0}2^{\sum\limits^4_1\!{}_s
f_s} {n!\over f_0!\prod\limits^4_1\!{}_sf_s!\prod\ov
R!c_{0000}!}c_{1200}!
{r_{2000}+r_{2100}\choose 
c_{2000},c_{2100}}
c_{2100}!
{r_{3000}+r_{3010}\choose 
c_{3000},c_{3010}}
c_{3010}!\\
\cdot{r_{1000}+r_{1200}\choose  c_{1000},c_{1200}}
{r_{2100}+r_{0100}\choose 
c_{2100},c_{0100}}
{r_{0010}+r_{3010}\choose 
c_{0010},c_{3010}}
{r_{0200}+r_{1200}\choose 
c_{0200},c_{1200}}\\
\cdot(f_1+r_{0000}+r_{3000}+r_{2000}+r_{0200})! 
(f_2+r_{0000}+r_{3000}+r_{0100}+r_{1000})!\\ 
\cdot(f_3+r_{0000}+r_{2000}+r_{1000}+r_{0010})! 
(f_4+r_{0000}+r_{0100}+r_{0010}+r_{0200})!
\endmultline
$$
that is the result already obtained by Pranesachar and others in 
[1]. If instead we apply (6.5) ,we obtain:
$$
\multline
K(4,n)=\sum_{f_0+f+\sum{\ov R}=n}\ \sum_{\ov
C-\ov D}(-1)^{f_0+r_{1000}+r_{0100}+r_{0010}+r_{3000}
+r_{2000}+r_{0200}}\\
6^{f_0}2^{f}
{n!\over f_0!\prod\ov R!\prod (\ov
C-\ov D)!}(r_{0000}+r_{3000}+r_{2000}+r_{0200})!\\
\cdot(r_{0000}+r_{3000}+r_{0100}+r_{1000})!
(r_{0000}+r_{2000}+r_{1000}+r_{0010})!\\
\cdot(r_{0000}+r_{0100} +r_{0010}+r_{0200})!
(r_{1000}+r_{1200})!(r_{2100}+r_{0100})!\\
\cdot(r_{0010}+r_{3010})!
(r_{0200}+r_{1200})!
(r_{2000}+r_{2100})!
(r_{3000}+r_{3010})!\\
\cdot{f+4r_{0000}+2\left(
r_{3000}+r_{2000}+r_{0100}+r_{1000}+r_{0200}+
r_{0010}\right)+3\choose f}\ .
\endmultline \tag7.4
$$
Which is an improvement on the results known up to now, since it
needs only 15 independent variables (the ten of $\ov R$, the three
of  $\ov C-\ov D$ and the two $f,f_0$) as compared with
the 18 of the formula of Pranesachar, Athreya and Singhi.

\head 8. Another point of view\endhead

In conclusion we want to show how the Theorem 5.1 can have
another interpretation which sheds light on its combinatory nature
 in a more profound way.

The circumstance --- which won't have escaped a careful reader
--- that the elements of  $R$ and of  $C$  are as many as those of
$\goto_m$, and that is  $m!$, is not casual. In fact if we  interpret
the indices  $\alpha_m,\dots, \alpha_1$ and $\beta_m,\dots,
\beta_1$ as the inversion tables of one of the permutations of 
$[m]$,  putting $b_i=\alpha_{m+1-i}$ (or $b_i=\beta_{m+1-i}$),
we will have two bijective maps between $C$ and $R$
and $\goto_m$, since $0\le \alpha_{m+1-i}\le m-i$.

Furthermore $\zeta(\aj)$ will be the same as the number of cycles of
$\sigma\in \goto_m$ which corresponds in this way to $r_\aj$.
However it is not true  --- as could be thought --- that the
components  $Z_l$ of $\aj$ correspond, in some way, to the cycles of
the permutation $\sigma$  corresponding to $r_\aj$.

To achieve this result we must introduce a new concept. Let us
take a  $\sigma\in\goto_m$ written in its standard representation
and put, $\forall \,i\in [m]\,k_i$ equal to $n+1-t$ where $t$ is the
element furthest on the right among those to the left of $i$
satisfying $t>i$ (or if  $i=s_h$, $k_i=m+1-s_t$ with $s_t>s_h$ and $t$
maximum); moreover we set $k_i=0$ if there are no
elements greater than $i$ on the left of  $i$. We say that 
$K(\sigma)=(k_1,\dots,k_m)$ is the ``covering table" of  $\sigma$. It
can be proved that the function $K(\sigma)$  is a bijection.
Furthermore it is clear that $0\le k_i\le m-i,\ \forall\, i\in [m]$, and
that, if $k_i=0$, $i$ is a  left-to-right maximum of the standard
representation of $\sigma$.

Now, if  we  put $k_i=\alpha_{m+1-j}$, we have that, $\forall\, i\in[m]$,
$0\le k_i\le m-i$ and therefore that $(k_1,\dots,k_m)$ can be
interpreted as the covering table of a  $S(r_\aj)\in\goto_m$. It can
be easily proved that  $S(r_\aj)$ is a bijection between $R$
and $\goto_m$ and that, in this case too, $\zeta(\aj)$ is the number
of the cycles of $S(r_\aj)$. Here however, if $s_hs_{h+1}\dots
s_{h+p}$ are the elements of a cycle of  $\sigma$ written in its
standard representation and if we take $k_{s_h},k_{s_{h+1}},\dots,
k_{s_{h+p}}$, we have that $\{\alpha_{m+1-s_h},
\alpha_{m+1-s_{h+1}},\dots,\alpha_{m+1-s_{h+p}}\}$ constitute a
component $Z_{m+1-s_h}$ of $\aj$.

In the light of this new bijective map, the results
obtained previously can be expressed in a new combinatory
language. In fact we can now consider the new variables $r_\sigma$
and $c_\vartheta$, the indices of which consist of elements of 
$\goto_m\,(\sigma,\vartheta\in\goto_m)$  and again indicate their
sets with  $R$ and $C$. Furthermore, writing
$\gamma\,|\,\sigma\in\goto_m$ to say that $\gamma$ is a cycle of 
$\sigma$, we can put $Q(\gamma)=\{r_\sigma\,|\,
\gamma\,|\,\sigma\}$ and $ q(\gamma)=\sum Q(\gamma)$;
corresponding meaning, going  from $r_\sigma$ to  $c_\vartheta$,
will have $\tilde Q(\gamma)$ and $\tilde q(\gamma)$. In this way
(5.6) can be reformulated like this: {\baselineskip 13pt
$$
\align
K(m,n)&=\sum_{\sum R=n}\hskip-.3truecm{}_{{}_R}\sum_{\sum \tilde
Q(\gamma)=q(\gamma)}\hskip-.7truecm{}_{{}_C}\ \
(-1)^{{}_{nm+\sum W'}}{\prod_{\gamma|\sigma}q(\gamma)!\over
\prod R!\prod C!}\tag8.1\\
 &=(-1)^{n(m-1)}\sum_{\sum R=n}\hskip-.2truecm{}_{{}_R}
\sum_{\sum \tilde
Q(\gamma)=q(\gamma)}\hskip-.7truecm{}_{{}_C}\ \
{\prod_{\gamma|\sigma}\left(- q(\gamma)\right)\t{\rm !` }\over
\prod R!\prod C!}
\endalign
$$
 where $W'$ indicates the set of all the $r_\sigma$
in which
$\sigma$ has an odd number of cycles.

The simplifications of section 6, which conduct us to (6.5),  can also
be read more clearly now. In fact $F_0$ consists of all the
$r_\sigma$ in which $\sigma$ is made up of only one cycle of order
$m$, while $F_s$, with $s\in[m]$,  is formed of those   $r_\sigma$ in
which $\sigma$ has a fixed point, made up of the element  $m+1-s$,
and a  cycle of order  $m-1$ which permutes the other elements
of  $[m]$.

Furthemore (8.1) reminds somehow the result attained by Gessel [3].

\head 9. The latin squares\endhead

When $m=n$,  we find ourselves facing the Latin
squares, much more famous than the Latin rectangles for their
applications in various branches of mathematics.

The number of  $n\times n$ Latin squares is usually indicated by
$L(n)$.

If we put $m=n$ in (5.6) and in (8.1), and, abandoning the
condition that the first row is in standard form, we multiply
everything by $n!$, we obtain the following elegant result which
allows us to  count of the number of Latin squares of any
order.

\proclaim{9.1 Theorem}
$$
L(n)=n!\sum_{\sum R=n}\!{}_{{}_R}\sum_{\tilde
Q(Z_l)=q(Z_l)}\hskip-.7truecm{}_{{}_C}\ \ {\prod\limits^n_0\!{}_l\prod\!{}_{Z_l}\left(-q(Z_l)
\right)\t{\rm !`}
\over\prod R!\prod C!}
=n!\sum_{\sum R=n}\hskip-.2truecm{}_{{}_R}\sum_{\tilde
Q(\gamma)=q(\gamma)}\hskip-.5truecm{}_{{}_C}\ \ {\prod\!{}_{\gamma|\sigma}
\left(-q(\gamma)\right)\t{\rm !`}
\over\prod R!\prod C!}\ .
\tag9.1
$$\endproclaim

\noindent In which  the  $2n!$ parameters  $R$ and $C$, and
their totals $q(Z_l)$ and $q(\gamma)$  previously
defined, appear.} 

\head References\endhead

\item{[1]} K.B. Athreya, C.R. Pranesachar and N.M. Singhi, On
the number of Latin rectangles and chromatic polynomials of
$L(K_{r,s})$, {\it European J. Combin.\/} {\bf 1} (1980), 9-17.

\item{[2]} P. Erdo\"s and I. Kaplansky, The asymptotic number of
Latin rectangles, {\it Amer. J. Math.\/} {\bf 68} (1946),
230-236.

\item{[3]} I.M. Gessel, Counting Latin rectangles, {\it Bull.
Amer. Math. Soc.\/} {\bf 16} (1987), 79-82.

\item{[4]} C.D. Godsil and B.D. McKay, Asymptotic enumeration
of Latin rectangles, {\it J. Comb. Th.\/} {\bf 48} (1990),
19-44.

\item{[5]} S.M. Kerawala, The enumeration of the Latin rectangle
of depth three by means of difference equations, {\it Bull. Calcutta
Math. Soc.\/} {\bf 33} (1941), 119-127.

\item{[6]} P.R. Montmort, Essai d'Analyse sur les Jeux de Hazard,
Paris (1708).

\item{[7]} W.O.J. Moser, The number of very reduced $4\times
n$ Latin rectangles, {\it Can. J. Math.\/} {\bf 19} (1967),
1011-1017.

\item{[8]} J.R. Nechvatal, Enumeration
of Latin rectangles, dissertation, University of Southern
California (1979).

\item{[9]} J.R. Nechvatal, Asymptotic enumeration
of generalized Latin rectangles, {\it Utilitas Math.\/} {\bf 20}
(1981), 273-292.

\item{[10]} R. Pranesachar, Enumeration
of Latin rectangles via SDR's, Combinatorics and Graph
Theory, (S.B. Rao, ed), {\it Lecture Notes in Math.\/} {\bf
885} (1981), 380-390.

\item{[11]} J. Riordan, Three-line  Latin rectangles II, {\it Amer.
Math. Month.\/} {\bf 53} (1946), 18-20.

\item{[12]} J. Riordan, Discordant permutations, {\it Scripta
Math. J. \/} {\bf 20} (1954), 14-23.

\item{[13]} R.P. Stanley, {\it Enumerative Combinatorics\/}, Vol.
I, Wadsworth and Brook\-s/Cole, (1986).

\item{[14]} C.M. Stein, Asymptotic evaluation of the number of
Latin rectangles, {\it J. Comb. Theory \/} {\bf (A) 25} (1978),
38-49.

\item{[15]} J. Touchard, Sur let probl\'eme de permutations, {\it
C. R. Acad. Sci. Paris \/} {\bf 198} (1934), 631-633.

\item{[16]} E.G. Whitehead Jr., Four-discordant permutations, {\it
J. Austral. Math. Soc. \/} {\bf A28} (1979), 369-377.

\item{[17]} K. Yamamoto, On the asymptotic  number of Latin
rectangles, {\it Jap. J.  Math. \/} {\bf 21} (1951),
113-119.

\item{[18]} K. Yamamoto, Structure polynomial of Latin
rectangles and its application to a combinatorial problem,
{\it Mem. Fac. Sci. Kyushu Univ.\/} {\bf A10} (1956), 1-13.

\bye